\documentclass[11pt]{amsart}
\usepackage{amsmath,amssymb}

\newtheorem{thm}{Theorem}[section]
\newtheorem{prop}[thm]{Proposition}
\newtheorem{lem}[thm]{Lemma}
\newtheorem{cor}[thm]{Corollary}

\newtheorem{question}[thm]{Question}
\newtheorem{rem}[thm]{Remark}
\newtheorem{defn}[thm]{Definition}
\newtheorem{exmp}[thm]{Example}

\theoremstyle{definition}
\newtheorem{dis}[thm]{Discussion}

\DeclareMathOperator{\U}{\mathcal U}

\DeclareMathOperator{\Spec}{Spec}
\DeclareMathOperator{\Max}{Max}

\DeclareMathOperator\Rad{Rad}

\DeclareMathOperator{\Z}{\mathbb Z}
\DeclareMathOperator{\C}{\mathcal C}
\DeclareMathOperator{\F}{\mathcal F}
\DeclareMathOperator{\G}{\mathcal G}

\DeclareMathOperator{\M}{\mathcal M}
\DeclareMathOperator{\PP}{\mathcal P}
\DeclareMathOperator{\Set}{\mathcal S}

\title[Well-Centered Extensions]
{Well-centered overrings  \\ of an integral domain}

\author[W. Heinzer and M. Roitman]{}



\begin{document}
\baselineskip 17pt

\maketitle

\begin{center}
{\bf William Heinzer}\\
Department of Mathematics,\\
Purdue University, West Lafayette, IN 47907-1395\\
E-mail: heinzer@math.purdue.edu

and

{
\bf Moshe Roitman
\footnote
{Part of this work was prepared while M. Roitman enjoyed the
hospitality of Purdue University.}
}\\
Department of Mathematics,\\
University of Haifa, Mount Carmel, Haifa 31905, Israel\\
E-mail: mroitman@math.haifa.ac.il\\

\end{center}

\begin{abstract} Let $A$ be an integral domain with field of fractions
$K$. We investigate the structure of the overrings $B \subseteq K$ of
$A$ that are well-centered on $A$ in the sense that each principal ideal
of $B$ is generated by an element of $A$. We consider the relation of
well-centeredness to the properties of flatness, localization and
sublocalization for $B$ over $A$. If $B = A[b]$ is a simple extension of
$A$, we prove that $B$ is a localization of $A$ if and only if $B$ is
flat and well-centered over $A$. If the integral closure of $A$ is a
Krull domain, in particular, if $A$ is Noetherian, we prove that every
finitely generated flat well-centered overring of $A$ is a localization
of $A$. We present examples of (non-finitely generated) flat
well-centered overrings of a Dedekind domain that are not localizations.
 \end{abstract}

\baselineskip 17pt

\section{Introduction.}

All rings we consider here are assumed to be commutative with unity.
If $R$ is a ring, we denote by $\U(R)$ the multiplicative group of
units of $R$.  If $A$ is an integral domain with field of
fractions $K$, we refer to  a
subring $B$  of $K$ with $A \subseteq B$ as an {\it overring } of $A$.

{\it Fix an integral domain $A$ with field of fractions 
$K$ and an  overring $B$ of $A$. }
\smallskip

  We say that $B$ is {\it well-centered on } $A$ if for each
$b \in B$ there exists a unit $u\in B$
such that $ub = a \in A$.  Thus, $B$ is well-centered on $A$ iff
each element of $B$ is an associate in $B$ of an
element of $A$ iff each
principal ideal of $B$ is generated by an element of $A$.

The overring  $B$ of $A$ is a {\it localization } of $A$ if $B = S^{-1}A = A_S$,
where $S$ is a multiplicatively closed subset of nonzero elements of $A$.
Thus $B$ is a localization of $A$ iff $B = A_{\U(B) \cap
A}$. A localization of $A$ is both flat over $A$ and well-centered on
$A$. Conversely, we prove in Theorem \ref{simple} that a simple flat
well-centered overring of an integral domain $A$ is a localization of
$A$. If the integral closure of $A$ is the intersection of a family of
valuation domains of finite character, we prove in Theorem
\ref{intclkrull} that every finitely generated flat well-centered
overring of $A$ is a localization of $A$. Thus every finitely generated
flat well-centered overring of an integral domain $A$ which is either
Krull or Noetherian is a localization of $A$ (Corollary \ref{fgnoe}). On
the other hand, we establish in Theorem \ref{exdedekind} the existence
of non-finitely generated flat well-centered overrings of a Dedekind
domain that are not localizations.

The overring $B$ of $A$ is  a {\it sublocalization } of $A$
if $B$ is an intersection of localizations of $A$. Thus $B$ is a
sublocalization of $A$ if and only if there exists a family
$\{S_{\lambda} \}_{\lambda \in \Lambda}$ of multiplicatively closed
subsets of nonzero elements of $A$ such that $B = \bigcap_{\lambda \in
\Lambda}A_{S_\lambda}$. It is well-known \cite{R}, \cite{GH1}
that a sublocalization $B$ of $A$ is an
intersection of localizations of $A$  at prime 
ideals. Indeed
 $\bigcap_{\lambda \in \Lambda}A_{S_\lambda} = \bigcap \{A_P : P \in \Spec
A \text{ and } P \cap S_{\lambda} = \emptyset \text{ for some } \lambda
\in \Lambda \}$ (see Discussion (\ref{inter})). 

A sublocalization $B$ of $A$ need be neither well-centered on $A$ nor flat
over $A$. We discuss in \S 2 the sublocalization condition in relation
to the properties of flatness and well-centeredness for an overring $B$
of $A$. We give in Corollary \ref{eachnoe} necessary and sufficient
conditions for each sublocalization overring of a Noetherian domain $A$
to be a localization of $A$.

We prove in Theorem \ref{fingenintcl} that every finitely generated
well-centered overring of an integrally closed domain is flat and
therefore, in particular, a sublocalization. In Example
\ref{wcfactorial} we establish the existence of a non-archimedean
well-centered overring of a factorial domain.

Our interest in the well-centered property of an overring of an integral
domain $A$ arose from conversations that the first author had with Jack
Ohm a number of years ago. The property arises naturally in relation to
results established by Ohm in Theorem 5.1 and Example 5.3 of \cite{O}.
M. Griffin in \cite[page 76]{Gr} defines well-centeredness
of a valuation $v$ with ring $B$ containing the 
domain $A$  in a manner equivalent to 
the definition of $B$ being well-centered
on $A$ given above. We thank Muhammad Zafrullah for 
pointing out to us this reference to Griffin. We also
thank the referee for several helpful suggestions that
have improved the paper.

\section{When a sublocalization is flat or a localization.}

Interesting work on the structure of flat overrings of an integral
domain has been done by Richman in \cite{R} and Akiba in \cite{A}.
Richman observes that an overring $B$ of $A$ is a flat $A$-module if
and only if $B_M = A_{M \cap A}$ for every maximal (or equivalently
prime) ideal $M$ of $B$ \cite[Theorem 2]{R}.  In particular, if $B$ is
a flat overring of $A$ then $B$ is a sublocalization of $A$.  The
converse of this result, however, is not true in general.  We indicate
below methods for obtaining sublocalizations $B$ of $A$ that fail to
be flat over $A$.

\begin{dis}\label{inter}

\
(1) If $B$ is a flat overring of $A$, then every ideal $J$ of $B$ is
extended from $A$. Indeed, for each maximal ideal $M$ of $B$ we have
$B_M=A_{M \cap A}$, hence $JB_M=JA_{M \cap A}
=(J\cap A)A_{M \cap A}=((J\cap A)B)B_M$.  Thus $J
= (J \cap A)B$. It is not true, however, that a flat overring $B$ of
an integral domain $A$ need be well-centered on $A$ (cf. Proposition
\ref{nontor} and Example \ref{exsimple}). The distinction is that
principal ideals of a flat overring $B$ need not be the extension of
principal ideals of $A$.

(2) If $S$ is a multiplicatively closed subset of an integral domain
$A$ with $0 \not\in S$, then
$$A_S = \bigcap \{ A_P : P \in \Spec A \text{ and } P \cap S =
\emptyset \}.
$$
Therefore if $\{S_{\lambda} \}_{\lambda \in \Lambda}$ is a family of
multiplicatively closed sets of nonzero elements of $A$ and $B =
\bigcap_{\lambda \in \Lambda}A_{S_\lambda}$, then
$$
B = \bigcap \{A_P : P \in \Spec A \text{ and } P \cap S_{\lambda} =
\emptyset \text{ for some } \lambda \in \Lambda \}.
$$
Thus  $B$ is a sublocalization over  $A$ if and only if
$$
B = \bigcap \{ A_P : P \in \Spec A \text{ and } B \subseteq A_P \}.
$$
In contrast with this characterization of a sublocalization, the
condition for each $P \in \Spec A$ that either $PB = B$ or $B
\subseteq A_P$ is, in general, stronger than the sublocalization
property.  Indeed, by \cite[Theorem 1]{R}, this latter property is
equivalent to flatness of $B$ over $A$.  Thus every flat overring is a
sublocalization. Hence every flat overring of an integrally closed
domain is again integrally closed \cite[Corollary, page 797]{R}.  Also
from Richman's characterization that $B$ is a flat overring of $A$ iff
for each $Q \in \Spec B$, we have $B_Q = A_{Q \cap A}$ \cite[Theorem
2]{R}, it follows that if $B$ is a quasilocal flat overring of $A$,
then $B$ is a localization of $A$.

(3) A useful observation is that if an overring $B \subseteq K$ of $A$
has one of the properties of being flat, well-centered, a
localization, or a sublocalization over $A$, then for each subring $C$
of $B$ with $A \subseteq C$, it follows that $B$ as an extension of
$C$ is, respectively, flat, well-centered, a localization, or a
sublocalization. This is easily seen in each case.

(4) If $B$ is a flat overring of $A$ and $C$ is a subring of $B$ with
$A \subseteq C$ such that $B$ is integral over $C$, then $C = B$. For
in this case $B$ is a flat integral overring of $C$, so by
\cite[Prop. 2]{R}, $C = B$.

(5) The localization, well-centered and flatness properties are
transitive in the sense that if $B$ is an overring of $A$ and $C$ is
an overring of $B$, then one of these properties holding for $B$ over
$A$ and for $C$ over $B$ implies the property also holds for $C$ over
$A$.

(6) The localization and flatness properties also behave well with
respect to compositum in the following sense: for an arbitrary
overring $C \subseteq K$ of $A$, if $B$ is a localization or a flat
overring of $A$, then $C[B]$ is, respectively, a localization or a
flat overring of $C$. For if $B = S^{-1}A$, then $C[B] = S^{-1}C$,
while for flatness if $Q \in \Spec C[B]$ and $P = Q \cap B$, then $B_P
= A_{P \cap A}$ implies $C[B]_Q = C_{Q \cap C}$.
\end{dis}

It would be interesting to know precise conditions for a Noetherian
integral domain to admit a non-Noetherian sublocalization overring.
In Corollary \ref{eachnoe}, we describe the class of Noetherian
domains $A$ for which each sublocalization over $A$ is a localization
of $A$. In particular, a Noetherian domain  in this class does not
admit a non-Noetherian sublocalization overring.

We begin with more general considerations.
We use $\Rad I$ to denote
the radical of an ideal $I$.

\begin{dis}\label{nonNoe}  If $R$ is a  ring, we
define $P \in \Spec R$ to be an {\it associated prime } of
an ideal $I$ of $R$ if there exists $a \in R$ such that $P$ is a
minimal prime over $(I:_Ra) = \{ r \in R : ra \in I \}$ \cite[page
289]{B}, \cite[page 92]{L}, \cite{BH}.  An integral domain $A$ has the
representation
$$A = \bigcap \{ A_P : P \text{ is an associated prime of a principal
ideal of } A \}
$$
\cite[Prop. 4]{BH}.  Moreover, if each principal ideal of $A$ has only
finitely many associates primes, then by \cite[Prop. 4]{BH} for $S$ a
multiplicatively closed subset of $A$, we have
$$
A_S = \bigcap \{A_P : P \text{ is an associated prime of a principal
ideal and } P \cap S = \emptyset \}.
$$
\end{dis}

\begin{lem}\label{containedminimal}
Let $P$ be a prime ideal of an integral domain $A$. Then the following
three properties are equivalent:

\begin{enumerate}

\item
For each family $\mathcal Q$ of prime ideals of $A$, if
$P \subseteq \bigcup_{Q \in \mathcal Q}Q$, then
$P \subseteq Q$ for some
$Q \in \mathcal Q$.

\item
For each family $\mathcal Q$ of minimal primes over principal
ideals of $A$, if $P \subseteq \bigcup_{Q \in \mathcal Q}Q$, then
$P \subseteq Q$ for some
$Q \in \mathcal Q$.

\item
$P$ is the radical of a principal ideal.

\end{enumerate}

\end{lem}

\begin{proof}

$(1\implies (2)$ 

Obvious.

$(2)\implies (1)$ 

Let $P\subseteq\bigcup_{Q\in\mathcal Q}Q$, where
$\mathcal Q$ is a set of prime ideals. Thus $P$ is contained in the
union of the set $\mathcal M$ of all minimal primes over principal
ideals contained in one of the primes $Q\in\mathcal Q$. Hence $P$ is
contained in some prime in $\mathcal M$ which is contained in a prime
$Q\in\mathcal Q$.

$(1)\implies (3)$

Let $\mathcal Q$ be the set of prime ideals of $A$ that do not contain $P$.
Thus $P \nsubseteq
\bigcup_{Q \in \mathcal Q}Q$. Let $c$ be an element in $P\setminus
\bigcup_{Q \in \mathcal Q}Q$. Since $P$ and $Ac$ are contained in the
same prime ideals, it follows that $P=\Rad(Ac)$.

$(3)\implies (1)$

Assume that $P=\Rad(Ac)$ for some element $c\in A$.
Let  $\mathcal Q$ be a family of prime ideals of $A$
so that $P \subseteq \bigcup_{Q \in \mathcal Q}Q$. Thus $c\in Q$ for some
prime ideal $Q\in\mathcal Q$, which implies that $P\subseteq Q$.

\end{proof}

We generalize below the theorem for Dedekind domains stated on
page~257 of \cite{GG} (see \cite{GO}).

\begin{thm}\label{critsubloc}
Let $A$ be an integral domain with field of fractions $K$, and let
$\mathcal P$ be a set of prime ideals in $A$. Consider the
sublocalization $B=\bigcap_{P\in\mathcal P}A_P$.  The following are
equivalent:
\begin{enumerate}
\item
$B$ is a localization  of $A$.

\item
If $x\in K\setminus A$, and $(A:_Ax)\subseteq \bigcup_{P\in\mathcal
P}P$,
then $(A:_Ax)\subseteq P$ for some $P\in\mathcal P$.

\noindent
Moreover, if each principal ideal of $A$ has only finitely many
associated primes, then the following condition is equivalent to the
two conditions above:

\item
If $Q$ is an associated prime of a principal ideal such that
$Q\subseteq \bigcup_{P\in\mathcal P}P$, then $Q\subseteq P$ for some
$P\in\mathcal P$.
\end{enumerate}

\end{thm}

\begin{proof}
$(1) \Longrightarrow (2):\ $ Assume that $B=A_S$ for some
multiplicative subset $S$ of $A$. Let $x\in K$ such that
$(A:_Ax)\subseteq \bigcup_{P\in\mathcal P}P$, thus $(A:_Ax)\cap
S=\emptyset$, hence $x\notin A_S=B$. Thus there exists a prime
$P\in\mathcal P$ such that $x\notin A_P$. It follows that
$(A:_Ax)\subseteq P$.

$(2) \Longrightarrow (1):\ $ Let $\ S=A\setminus
(\bigcup_{P\in\mathcal P}P)$. We prove that $B=A_S$. If $s\in S$, then
$s$ is a unit in $A_P$ for all $P\in\mathcal P$, hence $s$ is a unit
in $B$. It follows that $A_S\subseteq B$. On the other hand let $b\in
B\setminus A$, thus $(A:_A b)\nsubseteq P$ for all $P\in\mathcal
P$. By assumption $(A:_A b)\nsubseteq \bigcup_{P\in\mathcal P}P)$,
that is, $(A:_A b)\cap S\ne\emptyset$. It follows that $b\in A_S$.

Assume now that each principal ideal of $A$ has only finitely
many associated primes.

$(2) \Longrightarrow (3):\ $ Since principal ideals in $A$ have
only finitely many associated primes,  an associated prime of a
principal ideal is of the form $\Rad(A:_Ax)$ for some $x\in K$
\cite[Prop. 3.5]{HO}.

$(3) \Longrightarrow (2):\ $ Let $x\in K$ such that $(A:_Ax)\subseteq
\bigcup_{P\in\mathcal P}P$. By assumption, there are only  finitely
many prime ideals $Q_1,\dots, Q_n$ minimal over $(A:_Ax)$.  If none of
the primes $Q_i$ is contained in $\bigcup_{P\in\mathcal P}P$, then
choose an element $t_i\in Q_i\setminus \bigcup_{P\in\mathcal P}P$ for
each $i$. Thus for some positive integer $m$, we have
$(\prod_{i=1}^nt_i)^m\notin \bigcup_{P\in\mathcal P}P$, a
contradiction. Hence at least one of the ideals $Q_i$ is contained in
$\bigcup_{P\in\mathcal P}P$, which implies that $(A:_Ax)$ is contained
in $\bigcup_{P\in\mathcal P}P$.
\end{proof}

\begin{thm}\label{eachsubloc}
Let $A$ be an integral domain with field of fractions $K$. Each
sublocalization over $A$ is a localization of $A$ if and only if for
each $x\in K\setminus A$, the ideal $\Rad(A:_Ax)$ is the radical of
a principal ideal.

Moreover, if each principal ideal of $A$ has only  finitely many
associated primes, then each sublocalization of $A$ is a localization
iff each associated prime of a principal ideal is the radical of a
principal ideal.
\end{thm}

\begin{proof}
If each ideal of the form $\Rad(A:_Ax)$ is the radical of a principal
ideal, then each sublocalization of $A$ is a localization by
Theorem \ref{critsubloc}.

Conversely, assume that each sublocalization of $A$ is a localization
of $A$.  Let $x\in K\setminus A$. By Theorem \ref{critsubloc},
$(A:_Ax)$ is not contained in the union of the prime ideals not
containing $(A:_Ax)$. Let $c$ be an element in $(A:_Ax)$ that does not
belong to this union. Thus $(A:_Ax)$ and $Ac$ are contained in the same
prime ideals, which implies that $\Rad(A:_Ax)=\Rad(Ac)$.

Assume now that each principal ideal of $A$ has only  finitely many
associated primes, and that each sublocalization of $A$ is a
localization.  Let $P$ be a prime associated with a principal ideal of
$A$.  By Theorem \ref{critsubloc}, $P$ is not contained in a union of
primes not containing $P$. Hence, by Lemma \ref{containedminimal}, $P$
is the radical of a principal ideal.

Conversely, if each principal ideal of $A$ has only  finitely many
associated primes and if each associated prime of a principal ideal is
the radical of a principal ideal, then each sublocalization of $A$ is
a localization by Theorem \ref{critsubloc}.
\end{proof}

We apply the above results to various  classes of integral domains.
In Corollary \ref{eachmori} we describe the class
of Mori domains and the class of semi-Krull domains
for which each sublocalization is a localization. In
Corollary \ref{eachnoe} we characterize the
Noetherian domains having this property.

We recall that $A$ is a {\it Mori domain} if $A$ satisfies the
ascending chain condition on integral divisorial ideals
\cite{Ba}. In particular, a Mori domain satisfies the ascending
chain condition on principal ideals (a.c.c.p.). Examples of Mori
domains include factorial and Krull domains as well of course as
Noetherian domains.  An integral domain $A$ is {\it semi-Krull}
\cite{Ml}, if $A = \bigcap_PA_P$, where $P$ ranges over the set of
height-one primes of $A$, this intersection has finite character, and
for each height-one prime $P$, every nonzero ideal of $A_P$ contains a
power of $PA_P$.

A nonzero prime ideal of a Mori domain or a semi-Krull domain is an
associated prime of a principal ideal iff it is a prime divisorial
ideal (see \cite[Theorem 3.2]{Ba} and \cite[Theorem 1.7]{BGR}).  Thus
by Discussion \ref{nonNoe}, if $A_S$ is a localization of a Mori
domain $A$ or a semi-Krull domain $A$, then $A_S =
\bigcap_{P\in\mathcal P}A_P$, where $\mathcal P$ is the set of prime
divisorial ideals $P \in \Spec A$ such that $P \cap S = \emptyset$.
Therefore if $B$ is a sublocalization over $A$, then $B$ has
the form $B = \bigcap_{P\in\mathcal P}A_P$, where $\mathcal P$ is a
set of prime divisorial ideals in $A$.

Theorem \ref{critsubloc} implies:

\begin{cor}\label{critmori}
Let $A$ be a Mori domain or a semi-Krull domain and let $\mathcal P$
be a set of prime ideals in $A$. Consider the sublocalization
$B=\bigcap_{P\in\mathcal P}A_P$.  The following are equivalent:

\begin{enumerate}

\item
$B$ is a localization  of $A$.

\item
If $Q$ is a prime divisorial ideal of $A$ and $Q\subseteq
\bigcup_{P\in\mathcal P}P$, then $Q\subseteq P$ for
some $P\in\mathcal P$.

\end{enumerate}
\end{cor}

Theorem \ref{eachsubloc} implies:

\begin{cor}\label{eachmori}
Let $A$ be a Mori domain or a semi-Krull domain.  Each sublocalization
over $A$ is a localization of $A$ if and only if each prime divisorial
ideal of $A$ is the radical of a principal ideal.
\end{cor}

\begin{cor}\label{eachnoe}
Let $A$ be a Noetherian integral domain.  Each sublocalization over
$A$ is a localization of $A$ if and only if each associated prime of a
principal ideal of $A$ is the radical of a principal ideal.
In particular, if $A$ has these equivalent properties,
then nonzero principal ideals of $A$ have no embedded
associated primes.
\end{cor}

A Krull domain has torsion divisor class group iff each prime
divisorial ideal (that is, prime ideal of height one) is the radical
of a principal ideal. Hence Corollary \ref{eachmori} implies:

\begin{cor}\label{eachkrull}
A Krull domain $A$ has torsion divisor class group if and
if every  sublocalization over $A$ is a localization of $A$.
\end{cor}

\begin{cor} \label{eachonedim}
Let $A$ be a one-dimensional integral domain.  If each maximal ideal
of $A$ is the radical of a principal ideal, then every sublocalization
over $A$ is a localization of $A$. The converse holds if $A$ has
Noetherian prime spectrum.
\end{cor}

\begin{proof}
A commutative  ring has Noetherian spectrum iff each prime ideal is
the radical of a finitely generated ideal \cite{OP}. Thus a
one-dimensional integral domain has Noetherian spectrum iff
each nonzero element
is contained in  only finitely many maximal ideals iff
principal ideals have only finitely many
associated primes.  Thus
Corollary \ref{eachonedim}
follows from Theorem \ref{eachsubloc}.
\end{proof}

\begin{question}\label{sublocalflat}
What (Noetherian) integral domains $A$ have the property that every
sublocalization extension is flat?
\end{question}

For a one-dimensional integral domain with Noetherian spectrum we give
in Theorem \ref{1dimnoe} a complete answer to Question
\ref{sublocalflat}.

\begin{thm}\label{1dimnoe}
Suppose $A$ is a one-dimensional integral domain with Noetherian
spectrum.  Then every sublocalization over $A$ is flat over $A$.
\end{thm}

\begin{proof}
Let $B$ be a sublocalization over $A$.  We may assume that $B
\subsetneq K$, where $K$ is the field of fractions of $A$.  By
Discussion
\ref{inter} (2), there exists a family $\{P_\alpha \}$ of prime ideals
of $A$ such that $B = \bigcap_\alpha A_{P_\alpha}$. Since $\dim A =
1$, we may assume that each $P_\alpha$ is a maximal ideal of $A$. Let
$Q_\alpha = P_\alpha A_{P_\alpha} \cap B$.  We have $B_{Q_\alpha} =
A_{P_\alpha}$ and $B = \bigcap_\alpha B_{Q_\alpha}$. Since $A$ has
Noetherian spectrum, the family $\{B_{Q_\alpha} \}$ has finite
character in the sense that a nonzero element of $B$ is a unit in all
but finitely many of the $B_{Q_\alpha}$. To prove that $B$ is flat
over $A$, we show for each maximal ideal $Q$ of $B$ that $B_Q = A_{Q
\cap A}$. Let $P = Q \cap A$ and let $S = A\setminus P$.  By
\cite[Lemma 1.1]{HO2} we have $S^{-1}B =
\bigcap_\alpha(S^{-1}B_{Q_\alpha})$.  Since $B_{Q_\alpha}$ is a
one-dimensional quasilocal domain, $S^{-1}B_{Q_\alpha}$ is either
$B_{Q_\alpha}$ if $S \cap Q_\alpha\ne\emptyset$ or $K$ otherwise.
Since $A_{P_\alpha} = B_{Q_\alpha}$, we see that $Q_\alpha$ is the
unique prime of $B$ lying over $P_\alpha$.  Thus if $Q \ne Q_\alpha$,
then $S \cap Q_\alpha$ is nonempty and $S^{-1}B_{Q_\alpha} = K$.  If
this were true for each $\alpha$, then $S^{-1}B = \bigcap_\alpha
S^{-1}B_{Q_\alpha} = K$, but clearly $S^{-1}B \subseteq B_Q$, a
contradiction. Hence $Q = Q_\alpha$ for some $\alpha$ and therefore
$A_P = B_Q$.
\end{proof}

\section{Properties of flat and well-centered overrings.}

Richman observes \cite[Theorem 3]{R} that a flat overring of a
Noetherian domain is Noetherian.  There exist Noetherian integral
domains with non-Noetherian sublocalizations that are ideal
transforms (\cite{EH} and \cite[Theorem 3.2]{EHKR}). If $B$ is a
non-Noetherian ideal transform of a Noetherian domain $A$, then $B$ is
not flat over $A$ by the result of Richman mentioned
above. Proposition \ref{wcnoeth} shows that $B$ with these properties
also fails to be well-centered on $A$.

\begin{prop}\label{wcnoeth}
A well-centered extension of a Noetherian domain is Noetherian.
\end{prop}

\begin{proof}
If $B$ is well-centered on $A$, then every ideal of
$B$ is the extension of an ideal of $A$. Thus if $A$
is Noetherian, then every ideal of $B$ is finitely
generated and $B$ is also Noetherian.
\end{proof}

We observe in Theorem \ref{fingenintcl} that a finitely generated
well-centered overring of an integrally closed domain is a flat
extension. In the proof of this result we use Proposition \ref{max}
which holds for arbitrary well-centered extension rings.

\begin{prop}\label{max}
Let $S$ be a well-centered extension ring of a ring $R$. If $M$
is a maximal ideal of $R$ such that $MS \ne S$, then $MS$ is a maximal
ideal of $S$.
\end{prop}

\begin{proof}
We have a natural embedding $R/M\hookrightarrow S/MS$. Moreover
the fact that $S$ is well-centered over $R$ implies that
$S/MS$ is well-centered over $R/M$. Since a well-centered
extension of a field is a field, $S/MS$ is a field and $MS$
is a maximal ideal of $S$.
\end{proof}

For an extension ring $S$ of a ring $R$, we consider the following
condition that is in general weaker than the well-centered property.

\begin{defn}
An extension ring $S$ of a ring $R$ is said to be {\it almost
well-centered on } $R$ if for each $s \in S$ there exists a positive
integer $n$ depending on $s$ and an element $u \in \U(S)$ such that
$us^n \in R$.
\end{defn}

The following remark concerning almost well-centered
extensions is clear.

\begin{rem}\label{extensionideals}
If $S$ is an almost well-centered extension ring of a ring  $R$,
then for each ideal $J$ of $S$ we have
$\Rad J = \Rad(J \cap R)S$.
\end{rem}

In view of Remark \ref{extensionideals}, we have the
following analogue of Proposition \ref{max}.

\begin{prop}\label{radmax}
Let $S$ be an almost well-centered extension ring of a ring
$R$. If $M$ is a maximal ideal of $R$ such that $MS \ne S$, then
$\Rad{MS}$ is a maximal ideal of $S$.
\end{prop}

\begin{thm} \label{fingenintcl}
If $B$ is a finitely generated almost well-centered overring of $A$
and if $A$ is integrally closed in $B$, then $B$ is flat over $A$.  In
particular, every finitely generated almost well-centered overring of
an integrally closed domain $A$ is flat over $A$.
\end{thm}

\begin{proof}
Let $Q$ be a maximal ideal of $B$ and let $P= Q \cap A$.  By
Proposition \ref{radmax}, $\Rad(PB) = Q$. The Peskine-Evans version of
Zariski's Main Theorem \cite{P}, \cite{E} implies there exists $s \in
A \setminus P$ such that $A_s = B_s$.  In particular, $A_P =
B_Q$. Thus $B$ is flat over $A$.
\end{proof}

\begin{prop} \label{simpleintcl}
If $B = A[u]$ is a simple overring of $A$, where $u$ is a unit of $B$,
and if $A$ is integrally closed in $B$, then $B$ is a localization of
$A$.
\end{prop}

\begin{proof}
Since $u^{-1} \in B$ it follows that $u^{-1}$ is integral over $A$
\cite[Theorem 15]{Ka}.  Thus $u^{-1} \in A$ and $B$ is a localization
of $A$.
\end{proof}

\begin{cor}\label{simple_almost}
A simple almost well-centered overring of an integrally closed domain
is a localization.
\end{cor}

\begin{proof}
Let $B = A[b]$ be a simple almost well-centered overring of an
integrally closed domain $A$.  By Theorem \ref{fingenintcl}, $B$ is
flat over $A$.  Since $B$ is almost well-centered over $A$, there
exist a positive integer $n$ and a unit $u \in \U(B)$ such that $ub^n
= a \in A$. Thus $B$ is a flat integral overring of $A[b^n] = A[u]$.
By Discussion \ref{inter} (4), $B = A[u]$ and $B$ is a localization of
$A$.
\end{proof}

Theorem \ref{fingenintcl} and Corollary \ref{simple_almost} may fail
if $A$ is not integrally closed.  We use Proposition
\ref{wcen} to show in Example \ref{inext} the existence of Noetherian
integral domains that admit simple proper well-centered integral
overrings.  Corollary \ref{eachnoe} shows that  in  an integral
domains having this property there are principal ideals with
embedded associated prime ideals.

\begin{prop}\label{wcen}
Let $B$ be an integral domain of the form $B = K + M$, where $K$
is a field and $M$ is a nonzero maximal ideal of $B$.  If $A$ is a
subring of $B$
such that $M \subset A$, then $B$ is
well-centered on $A$.
\end{prop}

\begin{proof}
Let $b \in B$.  Then $b = k + m$, where $k \in K$ and $m \in M$. If $k
= 0$, then $b \in A$. If $k \ne 0$, then $k$ is a unit of $B$ and $a :
= b/k = 1 + (m/k) \in A$. Hence $B$ is well-centered over $A$.
\end{proof}

\begin{exmp} \label{inext}
A simple well-centered integral (thus not flat) proper overring $B$ of
a Noetherian integral domain $A$ such that $B$ is a sublocalization of
$A$.  Moreover, each height-one prime of $A$ is the radical of a
principal ideal.
\end{exmp}

Let $E = F(c)$ be a simple proper finite algebraic field extension,
let $B$ be the localized polynomial ring $E[X,Y]_{(X,Y)}$, let $M =
(X,Y)B$, and let $A = F + M$.  Then $A$ is Noetherian and $B=A[c]$ is
a simple, proper integral extension of $A$.  Hence $B$ is not flat as
an $A$-module \cite[Prop.  2]{R}. Proposition \ref{wcen} implies that
$B$ is well-centered on $A$.

Since $B$ is factorial, $B$ is the intersection of the rings $B_Q$ as
$Q$ ranges over the nonzero principal prime ideals of $B$.  For such
$Q$ we have $Q\subsetneq M\subset A$, thus $B\subseteq A_Q$, so
$B_Q=A_Q$. It follows that $B$ is a sublocalization over $A$.  Since
$B$ is a unique factorization domain, each height-one prime of $B$ is
principal.  Since $M \subset A$, each height-one prime of $A$ is the
radical of a principal ideal.  \qed

The following example where $B$ is not well-centered on $A$
illustrates restrictions on generalizing Proposition \ref{wcen}.
The  original Example \ref{3.11},
as noticed by Jung-Chen Liu and her student Jing-Ping Tsai,
is wrong. We provide here a correct example for the same statement
(actually, a little improvement since in the new version $M$
is a maximal ideal of both $A$ and $B$).

\begin{exmp} \label{3.11}
Integral domains of the form $A=A_0+M\subseteq B=B_0+M$, where
$A_0,B_0$ are subrings of $A$ and $B$, respectively, and $M$ is a
maximal ideal of $B$ such that $B$ is not almost well-centered on $A$.
\end{exmp}

Let $X$ be an indeterminate over the field $\mathbb Q$ of rational
numbers and let  $B=B_0=\mathbb Q[X]$. 
Let $M=(X^2-2)B, A_0=\mathbb Q$ and $A=A_0+M$. Clearly, $M$ is
a maximal ideal of both $A$ and $B=B_0+M$. No power of $X+1$ is
associate in $B$ with an element of $A$; otherwise, since $\U(B)=\mathbb
Q\smallsetminus\{0\}$, we obtain that some power of $X+1$ is in $A$:
$(X+1)^n\in \mathbb Q+(X^2-2)\mathbb Q[X]$. Put $X\to \sqrt 2$ in order
to obtain the contradiction $(\sqrt 2+1)^n\in\mathbb Q$ (for the 
binomial expansion implies $(\sqrt 2+1)^n \notin \mathbb Q$). Thus $B$ is not
almost well-centered over $A$. \qed

\begin{dis}\label{units}
Let $B$ be an overring of an integral domain $A$ and let $S=\U(B)\cap
A$.  Then $B = B_S$ is a well-centered overring of $A_S$ if and only
if $B$ is a well-centered overring of $A$. Moreover,
$\U(A_S)=\U(B_S)\cap A_S$, and $B$ is a localization of $A$ if and
only if $A_S=B$.  Thus in considering the question of whether an
overring $B$ of an integral domain $A$ is a localization of $A$, by
passing from the ring $A$ to its localization $A_{\U(B) \cap A}$, we
may assume that $\U(B) \cap A = \U(A)$.  The localization question is
then reduced to the question of whether $A = B$.  In general, if $B$
is a well-centered overring of $A$ which properly contains $A$, then
$\U(A) \subsetneq \U(B)$.  For if $b \in B \setminus A$ and $u \in
\U(B)$ is such that $ub \in A$, then $u^{-1} \not\in A$
so $u \in \U(B) \setminus \U(A)$.
\end{dis}

If $A$ is a Dedekind domain, then every overring $B$ of $A$ is a flat
$A$-module, thus a sublocalization over $A$. Moreover, we have:

\begin{prop}\label{nontor}
Let $A$ be a Dedekind domain. The following conditions are equivalent:

\begin{enumerate}

\item
$A$ has torsion divisor class group.

\item
Every overring of $A$ is a localization of $A$.

\item
Every overring of $A$ is well-centered on $A$.

\item
$A$ has no proper simple overring with the same set of units.

\end{enumerate}
\end{prop}

\begin{proof}
$(1)\iff(2):$ By Corollary~\ref{eachkrull}, (2) holds if and only if
each maximal ideal of $A$ is the radical of a principal ideal, and
this is equivalent to (1).

It is clear that $(2) \implies (3)$ and $(3) \implies (4)$.
Thus it remains to show:

$(4)\ \implies (2):\ $ Assume that (2) does not hold. Then $A$
has a maximal ideal $P$ that is not the radical of a principal
ideal. We claim that $B = A[P^{-1}]$ is a simple flat overring of $A$
with $\U(B)=\U(A)$. Indeed, if $b \in P^{-1} \setminus A$, we have $B
= A[b]$ since both of these rings are equal to $\bigcap \{A_Q : Q \in
\Spec A \text{ and } Q \ne P \}$.
Suppose there exists an element $u \in \U(B)\setminus \U(A)$. Then
$u$ is not a unit in $A_P$, but either $u$ or $u^{-1}$ is in $A_P$.
We may assume that $u \in A_P$,
thus $u \in PA_P$.  Then  $u\in A$ and
$\Rad{uA} = P$, a contradiction.
\end{proof}

We show in Theorem \ref{prufer} that if $B$ is a finitely generated
overring of a Dedekind domain $A$, then $B$ is a localization of $A$
iff $B$ is well-centered on $A$ iff $B$ is almost well-centered on
$A$. However, for overrings of a Dedekind domain having nontorsion
class group, we present in Theorem \ref{exdedekind} examples of
well-centered overrings that are not localizations and examples of
almost well-centered overrings that are not well-centered.

If $A$ is a Dedekind domain, we denote its class group by $\C(A)$; if
$I$ is a nonzero fractional ideal of $A$, we denote the ideal class of
$I$ by $\C_A(I)$, and if $\mathcal P$ is a subset of $\Max A$, we
denote the set $\{\C_A(P)\,|\, P\in\PP\}$ by $\C_A(\mathcal P)$. The
complement of a subset $\PP$ of $\Max A$ is denoted by
$\PP^c$. We denote the submonoid generated by a subset $S$ of a monoid
by $\M(S)$, and the subgroup generated by a subset $S$ of a group by
$\G(S)$. Thus, if $S$ is a set of nonzero fractional ideals of a
Dedekind domain $A$ viewed as a subset of the ideal monoid of $A$, we
have $\M(C_A(S))=\C_A(\M(S))$.

We recall that if $A$ is a Dedekind domain, and $B$ is an overring of
$A$, then there exists a unique set of maximal ideals $\mathcal P$ in
$A$ such
that $B = \bigcap \{ A_P : P \in \mathcal P\}.$ The overring $B$ of
$A$ can also be described as the compositum of the overrings
$A[Q^{-1}]$ such that $Q \in \Max A \setminus \mathcal P$.  Thus for
each $Q \in \Max A$ we have $QB = B$ if and only if $Q \in \mathcal
P^c $.

\begin{prop}\label{princ_class}
Let $A$ be a Dedekind domain with field of fractions $K$
and let $B \subsetneq K$ be an overring of $A$, thus
$$
B = \bigcap \{ A_P : P \in \mathcal P\}.
$$
for a unique  subset $\mathcal P$ of $\Max A$.  Let $J$ be a nonzero
ideal
of $B$. Then $J=IB$  where $I$ is an ideal of $A$ belonging to
$\M(\mathcal P))$. Moreover, we have
\begin{enumerate}
\item
$J$ is a principal ideal of $B \iff \C_A(I)\in \mathcal
G(\C_A(\mathcal P^c))$.
\item
$J$ is an extension of a principal ideal of $A \iff
\C_A(I)\in-\M(C_A(\PP^c))$.
\end{enumerate}
\end{prop}

\begin{proof}
Part (1)  follows from \cite[Corollary 3]{C}. For part (2),
assume first that there exists a principal ideal $I_0$ of $A$ such that
$IB=I_0B$. Since $I\in\M(\PP)$, it follows that $I_0=II_1$, where
$I_1\in\M(\PP^c)$. Thus $\C_A(I)=-\C_A(I_1)\in -\M(\C_A(\PP^c))$.

Conversely, let $\C_A(I)\in-\M(\C_A(\mathcal P^c))$.  There exists an
ideal $I_1\in\M(P^c)$ such that $II_1$ is a principal ideal of
$A$. Also $J=(II_1)B$.
\end{proof}

Proposition \ref{princ_class} implies:

\begin{cor}\label{well_class}
Let $A$ be a Dedekind domain with field of fractions $K$
and let $B \subsetneq K$ be an overring of $A$, thus
$$
B = \bigcap \{ A_P : P \in \mathcal P\}.
$$
for a unique  subset $\mathcal P$ of $\Max A$.
Then
\begin{enumerate}

\item
$B$ is a well-centered extension of $A \iff$
$$
(\C_A\M(\mathcal P))\cap\mathcal G(\C_A(\mathcal P^c)) \subseteq
-\C_A(\mathcal M(\mathcal P^c)).
$$

\item
$B$ is an almost well-centered extension of $A \iff $ each element
of\\ $\\M(C_A(\mathcal P))\cap\G(\C_A(\mathcal P^c))$ has a positive
integer multiple in $-\M(\C_A(\PP^c))$.

\end{enumerate}
\end{cor}

\begin{thm}\label{exdedekind}

\
\begin{enumerate}

\item
There exists a Dedekind domain $A$ having a well-centered overring
that is not a localization.

\item
There exists a Dedekind domain $A$ having an almost well-centered
overring that is not well-centered.

\end{enumerate}

Moreover, in each case the domain $A$ can be chosen so that it has
exactly two almost well-centered overrings  that are not
localizations of $A$, these two overrings being also the unique almost
well-centered overrings $D$ of $A$ such that $\U(D)\cap A=\U(A)$.

\end{thm}

\begin{proof}
We will use the well known result of Claborn \cite{C} that every
Abelian group is the ideal class group of a Dedekind domain, along
with the fact that for a countably generated Abelian group $G$ and a
nonempty subset $S$ of $G$, there exists a Dedekind domain $A$ with
class group $G$ such that $S = \{ \C(P) : P \in \Max A \}$ if and only
if $S$ generates $G$ as a monoid \cite[Theorem 5]{GHS}.

Let $A$ be a Dedekind domain having ideal class group the infinite
cyclic group $\Z$.  Define
$$
B = \bigcap \{ A_Q : Q \in \Max A \text{ and } \C(Q) \le 0 \}.
$$
Since the set $\{ \C(P) : P \in \Max A \}$ generates $\Z$ as a monoid,
there exists $P \in \Max A$ with $\C(P) > 0$. Thus $B$ is a proper
overring of $A$. For a nonzero nonunit $a \in A$, if $aA = P_1^{e_1}
\cdots P_n^{e_n}$ is the factorization of the principal ideal $aA$ as
a product of maximal ideals, then $0 = e_1\C(P_1) + \cdots +
e_n\C(P_n)$.  Therefore $\C(P_i) \le 0$ for at least one of the $P_i$.
It follows that $A \setminus \U(A) = \bigcup \{ Q : Q \in \Max A
\text{ and } \C(Q) \le 0 \}$. Since the maximal ideals of $B$ lie over
the ideals $Q$ of $A$ with $\C(Q)\le0$, we see that $B \setminus \U(B)
= \bigcup \{ QB : Q \in \Max A \text{ and } \C(Q) \le 0 \}$, hence
$\U(B) \cap A = \U(A)$.

By Corollary \ref{well_class}, $B$ is almost well-centered on $A$:
indeed, since there exists $P\in \PP^c$ with $C_A(P)>0$, each element
of $\M(\C_A(\mathcal P))$ has a power in $-\M(\C_A(\PP^c))$.
Moreover, if there exists $P \in \Max A$ with $\C(P) = 1$, by
Corollary \ref{well_class}, $B$ is well-centered on $A$.

To obtain an example where $B$ is almost well-centered but not
well-centered on $A$ we argue as follows. 
By \cite[Theorem 8]{GHS}, there exists a Dedekind domain $A$ with
class group $\Z$ such that $\{ \C(P) : P \in \Max A \} = \{ -1, 2, 3
\}$.  The overring
$$
B = \bigcap \{ A_Q : Q \in \Max A \text{ and } \C(Q) \le 0 \}
$$
is a principal ideal domain, since the primes $P \in \Max A$ such that
$PB = B$ generate $\Z$ as a group. Hence for $Q \in \Max A$ with
$\C(Q) = -1$, we have $QB = bB$ is a principal ideal that is not
generated by an element of $A$.

Next we show that for each Dedekind domain $A$ with ideal class group
$\Z$ as constructed above, there are precisely two proper almost
well-centered overrings $D$ of $A$ such that $\U(D) \cap A = \U(A)$.
These are the overring $B$ as defined above and $C = \bigcap \{ A_P :
\C(P) \ge 0 \}$. A proof that $A \subsetneq C$, $C$ is almost
well-centered over $A$, and that $\U(C) \cap A = \U(A)$ is similar to
that given above to show $B$ has these properties. Moreover, if $D$ is
an overring of $A$ such that $\U(D) \cap A = \U(A)$, then either $D
\subseteq B$ or $D \subseteq C$. For otherwise, either there exists a
$Q \in \Max A$ with $\C(Q) = 0$ such that $QD = D$ or there exist $P,
Q \in \Max A$ with $\C(P) = r > 0, \C(Q) = -s < 0$ and $PD=QD =D$. In
the first case, $Q = aA$ is principal and $a \in \U(D) \cap A
\setminus \U(A)$. In the second case $P^sQ^r = aA$ is principal and
again $a \in \U(D) \cap A \setminus \U(A)$.

It remains to show that if $A \subsetneq D \subsetneq B$ or
$A\subsetneq D \subsetneq C$, then $D$ is not almost well-centered
over $A$. If $A \subsetneq D \subsetneq B$, then the ideal class group
of $D$ is a proper homomorphic image of $\Z$ and hence a finite cyclic
group, thus each nonzero ideal of $D$ has a power that is a principal
ideal.  Since $D \subsetneq B$, there exists $P \in \Max A$ with
$\C(P) < 0$ such that $PD \in \Max D$. By Proposition \ref{princ_class}
(2),
no power of $PD$ is an extension of a principal ideal of $A$.
Therefore $D$ is not well-centered on $A$. The proof that an overring
$D$ of $A$ with $A \subsetneq D \subsetneq C$ is not well-centered on
$A$ is the same. Thus $B$ and $C$ are the unique proper almost
well-centered overrings of $A$ such that every nonunit of $A$ remains
a nonunit in the overring.

If $A$ has no principal maximal ideals, then $B$ and $C$ as
defined in the previous paragraph  are the unique
almost well-centered overrings of $A$ that
are not localizations of $A$. For if $D$ is a proper well-centered
overring of $A$ distinct from $B$ and $C$, then there exists a
localization $E$ of $A$ such that $A \subsetneq E \subseteq D$.
Since $A$ has no principal maximal ideals, the ideal class
group of $E$ is a proper homomorphic image of $\Z$.
Therefore $E$ has finite class group and every overring of
$E$ is a localization of $E$. Thus $D$ is a localization of $E$
and $E$ is a localization of $A$, so $D$ is a localization of $A$.
\end{proof}

\begin{prop}
Let $A$ be a Dedekind domain such that each ideal class
in the class group $C(A)$ of $A$ contains
a maximal ideal. If\ $C(A)$ is torsionfree, then each overring
of $A$ is an intersection of
two principal ideal domains that are well-centered overrings  of $A$.
\end{prop}

\begin{proof} Let $B=\bigcap_{P\in\PP}P$ be an overring of $A$, where
$\PP$ is a set of maximal ideals of $A$. Since $C(A)$ is torsionfree it
can be linearly ordered. With respect to a fixed linear order $\ge$  on $C(A)$,
define $B^+=\bigcap_{\{P\in\PP\text{ and }C_A(P)\ge0\}}A_P$ and
$B^-=\bigcap_{\{P\in\PP\text{ and }C_A(P)\le0\}}A_P$. Then $B=B^+\cap
B^-$, the empty intersection being defined as the field of fractions of
$A$. Since each ideal class of $A$ contains a prime ideal, Proposition
\ref{princ_class} implies that $B^+$ and $B^-$ are well-centered over
$A$ and that each prime ideal of $B^+$ and $B^-$ is the extension of a
principal ideal of $A$. Thus $B^+$ and $B^-$ are principal ideal domains
that are well-centered overrings of $A$ with $B = B^+ \cap B^-$.
\end{proof}

In \S 4 we use the following well-known general result characterizing
flat overrings, see for example \cite[Theorem 1]{A}. 
The implication $(1)\Longrightarrow (2) $ in Proposition
\ref{critflat} holds without assuming that $B$ is an overring of $A$,
cf. \cite[Exer. 22, page~47]{B}.

\begin{prop}\label{critflat}
Assume that $B$ is an overring of an integral domain $A$ and that
$S \subseteq B$ is  such that $B=A[S]$.
Then the following conditions are equivalent:
\begin{enumerate}
\item
$B$ is a flat extension of $A$.

\item
For any element  $s\in S$ we have $(A:_A s)B=B$.
\end{enumerate}
\end{prop}

If $B$ is well-centered over $A$, then $B = A[\U(B)]$. Thus the
following corollary of Proposition \ref{critflat} is
immediate.

\begin{cor}\label{critwelc}
Assume that $B$ is a well-centered overring of the
integral domain $A$.  Then the following conditions
are equivalent:
\begin{enumerate}

\item
$B$ is a flat extension of $A$.

\item
For each  unit $u\in B$  we have $(A:_A u)B=B$.
\end{enumerate}
\end{cor}

We recall that an integral domain $B$ is said to be {\it
Archimedean } if for each nonunit $b \in B$ we have
$\bigcap_{n=1}^\infty b^nB = (0)$.

\begin{rem}
If $B$ is a localization of an Archimedean domain $A$ such that the
conductor of $B$ in $A$
is nonzero, then $B=A$.
\end{rem}

Indeed, suppose  $B=A_S$ and let  $0\ne a\in (A:_AB)$.
Then for each  $s\in S$ we have $a/s^n \in  A$ for all $n\ge1$. 
Since $A$ is Archimedean, it follows
that  $s$ is a unit in $A$. Hence $B=A$.

\begin{prop}\label{conductor}
Suppose  $B$ is an overring of a  Mori integral
domain $A$.  If  the
conductor of $B$ in $A$ is nonzero and $B$ is
flat over  $A$, then $A=B$.
\end{prop}

\begin{proof}
Since $A$ is Mori and $(A:_AB)\ne0$, there exists a finite subset $F$
of $B$ such that $(A:_A B)=(A:_AF)$.
Since $B$ is flat over $A$,  Proposition \ref{critflat} implies that
$(A:_AF)B=B$, hence $(A:_AB) = (A:_AB)B = B$. Therefore  $A=B$.
\end{proof}

\begin{exmp}
If $A$ is not Mori, the conclusion of Proposition
\ref{conductor} need not hold.
\end{exmp}

Indeed, let $k$ be a field and let $R = k[X, Y]$ be a polynomial ring
over $k$.  Then $B := R[1/Y]$ is a localization of
$A := R + XB$. The conductor of $B$ in $A$ contains $XB$ and hence is
nonzero.  Moreover, $A \subsetneq B$ since $ Y^{-1}\in B\setminus
A$. \qed

The following structural result is proved by Querr\'e in \cite{Q}.

\begin{prop} \em \cite{Q} \em
If $A$ is a Mori domain and $B$ is a sublocalization over $A$, then
$B$ is also Mori.  In particular, a flat overring of a Mori domain is
again a Mori domain.
\end{prop}

We observe in Proposition \ref{wcnoeth} that a well-centered
overring of a Noetherian domain is Noetherian.  Example
\ref{wcfactorial} shows that in general the Mori property is not
preserved by well-centered overrings.  Indeed,  Example
\ref{wcfactorial} establishes the 
existence of a polynomial ring  $A$ over a field and a 
well-centered overring $B$ of $A$ that is
not Archimedean. In particular,  $B$ fails to satisfy a.c.c.p. and 
therefore is not Mori.

\begin{exmp}\label{wcfactorial}
A well-centered overring of a factorial domain (even of a polynomial
ring over a field) is not necessarily Archimedean.
\end{exmp}

Let $k$ be a field and let $a, c$ be two
independent indeterminates over $k$. Define
$$
T_0 = k[a, c, \{\frac{a}{c^n} :  n \ge 1 \}].
$$
Proceeding inductively,  define integral domains $T_m$ for $m \ge 1$
as follows: let $V_m = \{ v_{m,t} : t \text{ is a nonzero nonunit in }
T_{m-1} \}$ be a set of independent indeterminates over $T_{m-1}$ and
define $T_m := T_{m-1}[\{ v_{m,t}, \frac{1}{v_{m,t}} : v_{m,t} \in V_m
\}]$. Thus $T_m$ is a domain extension of $T_{m-1}$ obtained by
adjoining the indeterminates in $V_m$ along with their inverses.  Let
$V = \bigcup_{m=1}^\infty V_m$ and define $W$ to be the union of the
set $\{ a, c \}$ with the set $\{ tv_{m,t} : v_{m,t} \in V \}]$. The
elements of $W$ are algebraically independent over $k$. Thus $A := k[
W ]$ is a polynomial ring over the field $k$.  Define $B :=
\bigcup_{m=1}^\infty T_m$.  Since $T_0$ is an overring of $k[a, c]$,
we see that $B$ is an overring of $A$.  Since every element of
$T_{m-1}$ is an associate in $T_m$ to an element of $A$, it follows
that $B$ is well-centered on $A$.  The domain $B$ is not Archimedean
since $a/c^n \in B$ for all positive integers $n$ although $a,
c \in B$ and $c$ is not a unit in $B$;  indeed, $c \not\in \U(T_0)$ 
and $T_0$ is a retract of $B$ under the retraction over $T_0$
that sends each $v \in V$ to 1.  \qed

\section{Finitely generated well-centered extensions.}

The structure of a simple flat extension $S = R[s] = R[X]/I$ of a
commutative ring $R$ is considered in \cite{R}, \cite{V1}, \cite{V2},
\cite{OR1}, \cite{OR2}. Richman in \cite[Prop. 3]{R} shows that if $A$
is an integrally closed domain and $B = A[a/b]$ is a simple flat
overring of $A$, then $(a,b)A$ is an invertible ideal of $A$. We observe
in Theorems \ref{unit} that a simple flat overring $B$ generated by a
unit of $B$ is a localization of $A$. It follows (Corollary \ref{assoc}
and Corollary \ref{simple}) that well-centered simple flat overrings are
localizations.

\begin{thm}\label{unit}
Let $A$ be an integral domain and let $B=A[u]$
be a simple flat overring of $A$, where $u$ is a unit of $B$.
There exists a positive integer $m$ such that $u^{-r}\in A$
for all integers $r\ge m$.  Thus  $B$ is a localization of $A$.
\end{thm}

\begin{proof}
Since $u^{-1}\in A[u]$, the element $u^{-1}$ is integral over $A$,
hence $A[u^{-1}]$ is a finitely generated $A$-module.
Since $B$ is a flat extension of $A$, we have
$
(A:_A A[u^{-1}])B =B.
$
Hence
there exist  $c_0,\dots, c_m\in (A:_A A[u^{-1}])$ with
$1=c_0+ c_1u + \dots +c_mu^m$.
Thus for each integer $r \ge m$ we have
$$
u^{-r} = c_0u^{-r} + c_1u^{-r+1} + \cdots + c_mu^{-r+m} \in A.
$$
In particular, $u^{-m},u^{-m-1}\in A$. This implies
that $B=A[u^{m+1}]$ is a localization of
$A$.  \end{proof}

\begin{cor}\label{assoc}
Let $B = A[b]$ be a simple flat overring of an integral domain
$A$. The following are equivalent.

\begin{enumerate}
\item
$B$ is a localization  of $A$.
\item
$B$ is well-centered on $A$.
\item
$B$ is almost well-centered on $A$.
\item
The element $b$ is associate in $B$ with an element of $A$.
\item
Some power of  $b$ is associate in $B$ with an element of $A$.
\end{enumerate}
\end{cor}

\begin{proof}
It is enough to prove $(5)\ \Longrightarrow\ (1)$. Assume for some
positive integer $n$ that $b^n =au$ with $a\in A$ and $u\in
\U(B)$. Then $b^n \in A[u]$ implies $B$ is a flat integral overring of
$A[u]$.  Therefore $B = A[u]$ \cite[Prop. 2]{R}.  Hence by Theorem
\ref{unit}, $B$ is a localization of $A$.
\end{proof}

As an immediate consequence of either Theorem \ref{unit} or Corollary
\ref{assoc} we have:

\begin{cor}\label{simple}
If $B$ is a simple flat
well-centered overring of an integral
domain $A$, then $B$  is a localization of $A$.  \end{cor}

We present several additional  corollaries of
Theorem \ref{unit} concerning  finitely
generated flat overrings.

\begin{cor} \label{intclinside}
Let $B$ be a finitely generated flat overring of an integral domain
$A$ and let $A'$ denote the integral closure of $A$ in $B$. If $B$ is
a localization of $A'$, then $B$ is a localization of $A$.
\end{cor}

\begin{proof}
Since $B$ is finitely generated over $A$, if $B$ is a localization of
$A'$, then $B = A'[u]$ where $u^{-1} \in A'$. It follows that $B =
A[u][A']$ is an integral flat overring of $A[u]$. Therefore $B =
A[u]$. By Theorem \ref{unit}, $B$ is a localization of $A$.
\end{proof}

\begin{thm} \label{prufer}
Let $A$ be a Pr\"ufer domain with Noetherian spectrum (for example, a
Dedekind domain), and let $B$ be a finitely generated overring of
$A$. The following are equivalent.

\begin{enumerate}
\item
$B$ is a localization  of $A$.

\item
$B$ is well-centered on $A$.
\item
$B$ is almost well-centered on $A$.
\end{enumerate}

\end{thm}

\begin{proof}
It is enough to prove $(3)\ \Longrightarrow\ (1)$.  Assume that $B$ is
almost well-centered on $A$.  By \cite[Corollary 5.6]{GH2}, $B = A[b]$
is a simple extension.  Since every overring of a Pr\"ufer domain is
flat, we obtain by Corollary \ref{assoc} that $B$ is a localization of
$A$.
\end{proof}

In Proposition \ref{nontor}, we present  examples of Dedekind
domains $A \subset B$ such that $B$ is a proper simple flat overring
of $A$ and $\U(A) = \U(B)$.  Example \ref{exsimple} provides a more
explicit construction of this type and also shows that the condition
that $u$ is a unit in $B$ is essential in Theorem~\ref{unit}.

\begin{exmp}\label{exsimple}
An example of a simple flat overring $B$ of an integrally closed
domain $A$ such that $A\subsetneq B$ and $\U(A)=\U(B)$.
\end{exmp}

Let $X,Y$ and $Z$ be indeterminates over a field $k$.  Set
\begin{equation*}
\begin{split}
A&=k[X,Y,XZ,YZ,\frac 1{X+YZ}],\\
B&=k[X,Y,Z,\frac 1{X+YZ}].
\end{split}\end{equation*}

Clearly $A$ and $B$ have the same field of fractions  $k(X,Y,Z)$ and
$B=A[Z]$. To see that $Z\in B\setminus A$, observe that the
$k[Z]$-algebra homomorphism defined by setting $Y = \frac{1}{Z}$ and $X
=
0$ maps $A$ to $k$.  Also
$\U(A)=\U(B)=\left\{ a(X+YZ)^m\,|\,a\in k\setminus\{0\},m\in\mathbb
Z\right\}$.  Since $A$ is a localization of the integrally closed
domain $k[X, Y, XZ, YZ]$, we see that $A$ is integrally closed.
Finally, $B$ is a flat extension of $A$ by Proposition \ref{critflat}
since $XZ, YZ \in A$ and the ideal $(X,Y)B = B$.

\begin{question}\label{fgwell}
Under what conditions on $A$ is every finitely generated
well-centered overring of $A$
a localization of $A$?
\end{question}

If $A$ is Noetherian, it follows from Corollary \ref{fgnoe} that every
finitely generated flat well-centered overring of $A$ is a
localization of $A$.  In a situation where Question \ref{fgwell} has a
positive answer, it follows that the finitely generated overring is
actually a simple extension, for if $B$ is a finitely generated
overring of $A$ that is a localization of $A$, then $B$ is a simple
extension of $A$.

\begin{rem}
Let $B=A[u, v]$ be a flat overring of a domain $A$, where $v\in\U(B)$.
Then $B=A[u,\dfrac 1 {f(u)}]$ for some polynomial $f(X)\in A[X]$.
\end{rem}

Indeed, $B$ is a localization of $A[u]$.

\begin{prop}\label{u1/u}
Let $B=A[u, \dfrac 1 u]$ be a flat overring of a domain $A$, where
$u\in\U(B)$.  Then $B=A[u+\dfrac 1 u]$ . Moreover, if $B$ is
well-centered over $A$, then $B$ is a localization of $A$.
\end{prop}

\begin{proof}
Let $C=A[u+\dfrac 1 u]$.  Since $B=C[u]=C[\dfrac 1 u]$,
we obtain by Theorem \ref{unit} that $u^{-n},u^n\in C$ for
sufficiently large $n$.  Hence $u,\dfrac 1 u\in C$, which implies
$C=B$.  By Corollary \ref{assoc}, if $B$ is well-centered over $A$,
then $B$ is a localization of $A$.
\end{proof}

We extend Proposition \ref{u1/u} as follows:

\begin{prop}
Let $A$ be an integral domain and let $B=A[u,\dfrac 1{f(u)}]$ be a
flat well-centered overring of $A$, where $f(X)$ is a monic polynomial
in $A[X]$, and $u,f(u)\in \U(B)$.  Then $B$ is a localization of $A$.
\end{prop}

\begin{proof}
Since $f$ is monic, $B$ is integral over 
$C:=A[f(u),\dfrac 1 {f(u)}]$. Thus $B$ is flat and integral
over $C$ and therefore $B = C$. Thus $B = C$ is flat and 
well-centered over $A$. Proposition \ref{u1/u} implies that
$B = A[f(u) + \dfrac 1 {f(u)}]$ and that $B$ is a localization
of $A$. 
\end{proof}

\begin{question}
Under what conditions on an integral domain $A$
is every flat overring of $A$ well-centered on $A$?
\end{question}

\begin{dis}\label{akiba}
Akiba in \cite{A} constructs an interesting example where $A$ is a
2-dimensional normal excellent local domain, $P$ is a height-one prime
of $A$ that is not the radical of a principal ideal, and $B =
\bigcup_{n=1}^\infty P^{-n}$ is the ideal transform of $A$ at
$P$. Thus $B = \bigcap_QA_Q$, where the intersection ranges over all
the height-one primes of $A$ other than $P$.  Akiba proves that $PB =
B$. It follows that $B$ is flat and finitely generated over $A$, but
not a localization of $A$.

We observe that $B$ is not almost well-centered over $A$.  Indeed,
assume that $B$ is almost well-centered over $A$ , and let $b \in B
\setminus A$. Thus $ub^m \in A$ for some unit $u$ of $B$ and $m\ge1$.
Hence $u \in \U(A_Q)$ for each height-one prime $Q \ne P$ of $A$.
Since $A$ is normal, we have $b^m \in B\setminus A$, thus $b^m \not\in
A_P$. It follows that $u \in PA_P$. Therefore $u \in A$ and $\sqrt{uA}
= P$. This contradicts the fact that $P$ is not the radical of a
principal ideal. We conclude that $B$ is not almost well-centered on
$A$.

We observe that $B$ is not a simple extension of $A$. Moreover, for
every nonzero nonunit $b \in B$ we have $C := A + bB \subsetneq B$.
This follows because $PB = B$ implies $\dim B = 1$ and $\dim(B/bB) =
0$. However, $C/bB \cong A/(bB \cap A)$ and $\dim(A/(bB \cap A)) = 1$.
\end{dis}

\begin{thm} \label{repasfinint}
Let $B$ be a well-centered overring of an integral domain $A$. If
there exist finitely many valuation overrings $V_1, \dots, V_n$ of $A$
such that $A = B \cap V_1 \cap \cdots \cap V_n$, then $B$ is a
localization of $A$.
\end{thm}

\begin{proof}
For $S$ a multiplicatively closed subset of $A$,
we have
$$
S^{-1}A  = S^{-1}B \cap S^{-1}V_1  \cap \dots \cap S^{-1}V_n,
$$
so by replacing $A$ by its localization $(\U(B) \cap A)^{-1}A$, we may
assume that $\U(B) \cap A = \U(A)$.  If $B \subseteq V_i$, then $V_i$
may be deleted in the representation $A = B  \cap
(\bigcap_{i=1}^n V_i)$. Thus we may assume that 
$B \not\subseteq V_i$ for each $i$.  We
prove that after these reductions we have $A = B$, i.e., the set $\{
V_i \}$ is empty. Assume not, 
then for each $1\le i\le n$ choose $b_i \in B$
such that $b_i \not\in V_i$.  By \cite[Lemma 5.4]{GH2}, there exist
positive integers $e_1, \dots,e_n$ such that $b := b_1^{e_1} +
b_2^{e_2} + \cdots b_n^{e_n} \not\in V_i$, thus $b^{-1} \in V_i$ for
each $i = 1, \dots, n$. Since $B$ is well-centered over $A$, there
exists $u \in \U(B)$ such that $ub \in A$. Since $b \not\in V_i$, 
we have  $u\in V_i$ for all
$i$.  Therefore $u \in B \cap (\cap_{i=1}^nV_i = A$. It follows that 
$u \in A \cap \U(B) = \U(A)$ and $u^{-1} \in A$. Hence
$b\in A$, a contradiction.
\end{proof}

\begin{lem}\label{locint}
Let $B$ be a finitely generated flat overring of an
integral  domain $A$ and let $C$ be an integral
overring of $A$. The following conditions are equivalent.

\begin{enumerate}

\item
$B$ is a localization of $A$.

\item
$B$ is a localization of $C\cap B$.

\end{enumerate}
\end{lem}
\begin{proof}
Clearly $(1) \implies (2)$. Assume (2). Then 
$B = (B \cap C)[u]$, where $u^{-1} \in B \cap C$. 
Since $B \cap C$ is integral over $A$, it 
follows that $B$ is flat and integral over $A[u]$.
Therefore $B = A[u]$. 
By Theorem \ref{unit}, $B$ is a
localization of $A$.  
\end{proof}

\begin{thm} \label{intclkrull}
Let $A$ be an integral domain for which the integral closure $A'$
has a representation $A' = \bigcap_{V\in\mathcal V}V$,
where $\mathcal V$ is a family of valuation overrings of $A$
of finite character.  If $B$ is a
finitely generated flat well-centered overring of $A$, then $B$ is a
localization of $A$.
\end{thm}

\begin{proof}
Since $B$ is finitely generated over $A$, we have $B \subseteq V$ for
all but finitely many domains $V\in\mathcal V$. Let $V_1,\dots,V_n$ be
the domains in $\mathcal V$ that do not contain $B$. Thus $A'\cap
B=B\cap(\bigcap_{i=1}^nV_i)$. By Theorem \ref{repasfinint}, $B$ is a
localization of $A'\cap B$. By Lemma \ref{locint}, $B$ is a
localization of $A$.
\end{proof}

It is well known that the integral closure of a Noetherian domain is a
Krull domain \cite[(33.10)]{N}. Therefore Corollary \ref{fgnoe} is an
immediate consequence of Theorem \ref{intclkrull}.

\begin{cor}\label{fgnoe}
Let $A$ be an integral domain for which the integral closure $A'$ is a
Krull domain. If $B$ is a finitely generated flat well-centered
overring of $A$, then $B$ is a localization of $A$. In particular, a
finitely generated flat well-centered overring of a Noetherian
integral domain $A$ is a localization of $A$.
\end{cor}

\begin{dis}\label{PB=B}
Let $A$ be an integral domain with field of fractions $K$.
Suppose $B=A[b_1,\dots,b_n]$ is a finitely generated overring of
$A$.  Let $I_j = (A:_A b_j)$ be the denominator ideal
of $b_j$ and let $I = \cap_{j=1}^n I_j$.  The overring
$C := \{ x \in K : xI^n \subseteq A \text{ for some integer } n \ge 1
\}$
is called the $I$-transform of $A$. This construction was  first
introduced by
Nagata \cite{N2} in his work on the 14-th problem of Hilbert. It is
clear that $B \subseteq C$ and that $C$ is the $IB$-transform of $B$.
Nagata observes \cite[Lemma 3, page 58]{N2} that there is a one-to-one
correspondence between the prime ideals $Q$ of $C$ not containing $I$
and the prime ideals $P$ of $A$ not containing $I$ effected by defining
$Q \cap A = P$. Moreover, it then follows that  $A_P = C_Q$. In
particular,
if $IB = B$, then $B = C$ is flat over $A$
and there is a one-to-one correspondence
between the prime ideals $Q$ of $B$ and the prime ideals $P$ of $A$
not containing $I$, the correspondence defined by $Q \cap A = P$. Thus
if
$B = A[b_1, \ldots, b_n]$ is a flat overring of $A$
and $P \in \Spec A$,
then the following are equivalent:
\begin{enumerate}

\item
$PB=B$.

\item
$P$ contains the  ideal $I_j = (A:_Ab_j)$ for some $j \in \{1, \ldots,
n\}$.

\item
$P$ contains the ideal $I = \cap_{j=1}^n I_j$.

\end{enumerate}
\end{dis}

\begin{thm}\label{finite1}
Let $B$ be a well-centered flat overring of an integral domain $A$.
If there exists a finite set $\F$ of height-one prime ideals of $A$
such that $A = B \cap \bigcap_{P\in \F}A_P$, then $B$ is a
localization of $A$.
\end{thm}

\begin{proof}
Since $P \in \F$ has height-one, for 
$S$ a multiplicatively closed subset of $A$ 
either $S^{-1}A_P = A_P$  or $S^{-1}A_P = K$, the field 
of fractions of $A$. Therefore 
$$ 
S^{-1}A = S^{-1}B \cap \bigcap_P(S^{-1}A_P),
$$ 
where the intersection is over all  $P\in\F$ such that
$P\cap S=\emptyset$.  By replacing $A$ by its localization $(\U(B)
\cap A)^{-1}A$, we may assume that $\U(B) \cap A = \U(A)$
and that $B \not\subseteq A_P$ for each $P \in \F$.
After this reduction, we claim that $A = B$, i.e., that 
$\F = \emptyset$. 
Suppose  $\F\ne\emptyset$. Since $B$ is flat
over $A$, for each $P\in \F$ we have $PB = B$.
Let $c$ be a nonzero element in $\bigcap_{P\in\F}P$ 
and consider the ring $B/cB$ and its subring 
$R = A/(cB \cap A)$. Since $PB = B$ and since
every minimal prime of the ring $R$ is the contraction
of a prime ideal of $B/cB$, we have
$cB\cap A\nsubseteq P$ for each 
$P\in\F$.  Thus there exists an element $s\in A\setminus
\bigcup_{P\in\F}P$, so that $s/c \in B$.  Since $B$ is
well-centered over $A$, there exists $u \in \U(B)$ such that 
$us/c=a \in A$.  Thus $u= ac/s \in A_P$ for all $P\in\F$.
Therefore $u \in A \cap \U(B) = \U(A)$.  Hence $s/c \in A$,
but $ s/c \notin A_P$, a contradiction.
\end{proof}

\begin{thm}\label{height1}
If each  nonzero principal ideal of  the integral domain
$A$ has only finitely many associated primes
and each of these associated primes is  of height $1$, 
then every  finitely generated flat well-centered
overring of $A$ is a localization of
$A$.  \end{thm}

\begin{proof}
Let $B=A[b_1,\dots,b_n]$ be a finitely generated flat well-centered
overring of  $A$.  To prove that $B$ is a localization of $A$, we
may assume that $\U(B)\cap A=\U(A)$, and then we have to show that
$A=B$.

Since $B$ is a sublocalization of $A$, by \cite[Prop.  4]{BH},
$B=\bigcap_{P\in \Set}A_P$, where $\Set$ is the set of prime ideals
$P$ of height $1$ of $A$ so that $PB\ne B$.  Let $\F$ be the set of
prime ideals of height $1$ in $A$ such that $PB=B$.  Then 
$A=B\cap\bigcap_{P\in\F}A_P$.  By Discussion \ref{PB=B}, the set $\F$
is finite.  Hence by Theorem \ref{finite1},  $B$ is a localization of $A$.
\end{proof}

\begin{prop}\label{multset}
Let $B$ be a well-centered overring of an integral domain $A$. If 
$S$ is  a multiplicative closed subset of $A$ such that $A=A_S\cap B$
and such that $B_S$ is a localization of $A_S$, then  $B$
is a localization of $A$. 
\end{prop}

\begin{proof}
Let $b\in B$. There exists
an element $t\in A_S\cap\U(B_S)$ such that $tb\in A_S$. We may assume
that $t\in A$, thus $tb\in A_S\cap B = A$. 
Since $t^{-1} \in B_S$, there exists $s \in S$ such that
$st^{-1} \in B$. Since $B$ is well-centered over $A$, there
exists $u \in \U(B)$ such that $ust^{-1} = a \in A$. 
Then $u = at/s \in A_S \cap B = A$ and 
$ub = atb/s \in A_S \cap B = A$. We have shown 
for each $b \in B$ there exists $u \in A \cap \U(B)$
such that $ub \in A$. Therefore $B$ is a localization of $A$. 
\end{proof}

\begin{cor}\label{multid}
Let $B$ be a well-centered overring of an integral domain $A$, let
$I$ be a proper ideal of $A$, and let $S = 1 + I$.
If for each $b\in B$ and $c \in I$
there exists an integer $n \ge 1$ such that $c^nb\in A$ and if  
$B_S$ is a localization of $A_S$, then $B$ is a localization of $A$.
\end{cor}

\begin{proof}
The corollary follows from Proposition \ref{multset} since $A=A_S\cap B$.
\end{proof}

\begin{thm}\label{onedim}
Every finitely
generated flat well-centered overring of a one-dimensional
integral domain $A$ is a localization of $A$.
\end{thm}

\begin{proof}
Let $B=A[b_1,\dots,b_n]$ be a finitely generated flat well-centered
overring of $A$ and let  $I = \bigcap_{j=1}^n (A:_A b_j)$. 
Then  $IB=B$ by flatness. Let
$S=1+I$. Then  $IA_S$ is contained in the Jacobson radical of $A_S$.
Since $\dim A_S \le 1$, $IA_S$ is contained in every nonzero 
prime ideal of $A_S$.  Since $IB=B$, it follows by 
\cite[Theorem 2 or 3]{R}, that $B_S$ is the field of fractions of $A_S$.
By Corollary \ref{multid}, $B$ is a localization of $A$.
\end{proof}

An interesting question that remains open is whether a finitely
generated flat well-centered 
overring of an integral domain $A$ is always a localization of $A$.

\end{document}